\pdfoutput=1
\RequirePackage{ifpdf}
\ifpdf 
\documentclass[pdftex]{sigma}
\else
\documentclass{sigma}
\fi

\usepackage{mathrsfs,bm,mathtools,accents,bbm}

\usepackage{pgf}
\usepackage{float}

\usepackage{xcolor}

\DeclareFontFamily{U}{stixextrai}{}
\DeclareFontShape{U}{stixextrai}{m}{n}
 { <-> stix-extra1 }{}

\newcommand*{\parallelogram}{%
 \rlap{\rotatebox{-30}{\rule[.05ex]{.4pt}{.77em}}}%
 \kern.04em%
 \rlap{\kern.36em\raisebox{0.649519052835em}{\rule{.6em}{.4pt}}}%
 \rule{.6em}{.4pt}\kern-.04em%
 \rotatebox{-30}{\rule[.05ex]{.4pt}{.77em}}}

\DeclareMathOperator*{\Res}{Res}

\DeclareMathOperator*{\Vol}{Vol}

\DeclareMathOperator*{\Arg}{Arg}

\DeclareMathOperator*{\sign}{sign}

\newcommand*{\HG}[4]{#1\left(\begin{array}{@{}c|}
#2 \\ #3
\end{array}\,\,
#4
\right)}

\usepackage{tikz}
\usetikzlibrary{tikzmark,arrows,matrix,decorations.pathreplacing,calc,decorations.markings}

\pgfkeys{tikz/mymatrixenv/.style={decoration=brace,every left delimiter/.style={xshift=3pt},every right delimiter/.style={xshift=-3pt}}}
\pgfkeys{tikz/mymatrix/.style={matrix of math nodes,left delimiter={(},right delimiter={)},inner sep=2pt,column sep=1em,row sep=0.5em,nodes={inner sep=0pt}}}
\pgfkeys{tikz/mymatrixbrace/.style={decorate,thick}}

\makeatletter
\newcommand*{\bigominus}{\DOTSB\bigominus@\slimits@}

\newcommand{\bigominus@}{\mathop{\mathpalette\bigominus@@\relax}}

\newcommand{\bigominus@@}[2]{%
 \vcenter{\hbox{%
 \sbox\z@{$\m@th#1\bigoplus$}%
 \resizebox{\wd\z@}{!}{$\m@th#1\bm{\ominus}$}%
 }}%
}
\makeatother

\DeclareFontFamily{U}{mathx}{\hyphenchar\font45}
\DeclareFontShape{U}{mathx}{m}{n}{
 <5> <6> <7> <8> <9> <10>
 <10.95> <12> <14.4> <17.28> <20.74> <24.88>
 mathx10
 }{}
\DeclareSymbolFont{mathx}{U}{mathx}{m}{n}
\DeclareFontSubstitution{U}{mathx}{m}{n}

\numberwithin{equation}{section}

\newtheorem{Theorem}{Theorem}[section]
\newtheorem*{Theorem*}{Theorem}
\newtheorem{Corollary}[Theorem]{Corollary}
\newtheorem{Lemma}[Theorem]{Lemma}
\newtheorem{Proposition}[Theorem]{Proposition}
 { \theoremstyle{definition}
\newtheorem{Definition}[Theorem]{Definition}
\newtheorem{Note}[Theorem]{Note}
\newtheorem{Example}[Theorem]{Example}
\newtheorem{Remark}[Theorem]{Remark} }

\begin{document}
\allowdisplaybreaks

\newcommand{\arXivNumber}{2011.00707}

\renewcommand{\PaperNumber}{048}

\FirstPageHeading

\ShortArticleName{On the Monodromy Invariant Hermitian Form for $A$-Hypergeometric Systems}

\ArticleName{On the Monodromy Invariant Hermitian Form\\ for $\boldsymbol{A}$-Hypergeometric Systems}

\Author{Carlo VERSCHOOR}

\AuthorNameForHeading{C.~Verschoor}

\Address{Department of Mathematics, Utrecht University, Utrecht,\\ Budapestlaan 6, 3580 TA, The Netherlands}
\Email{\href{mailto:carlovmm@gmail.com}{carlovmm@gmail.com}}

\ArticleDates{Received August 24, 2021, in final form June 22, 2022; Published online June 30, 2022}

\Abstract{We will give an explicit construction of the invariant Hermitian form for the monodromy of an $A$-hypergeometric system given that there is a~Mellin--Barnes basis of solutions.}


\Classification{14D05; 33C70}

\section{Introduction}

$A$-hypergeometric functions were introduced by Gelfand, Kapranov and Zelevinsky \cite{GKZ87, GKZ90, GKZ88, GKZ89} to give us a framework to study hypergeometric functions in more generality. Many classical hypergeometric functions can be described in terms of $A$-hypergeometric functions. Examples include Gauss' hypergeometric function ${}_2F_1$, the generalized hypergeometric functions of the type ${}_{n}F_{n-1}$, Appell's hypergeometric functions $F_1$, $F_2$, $F_3$, $F_4$ \cite{App80,App82}, Horn's hypergeometric functions~$G_1$, $G_2$, $G_3$, $H_1$, $H_2$, $H_3$, $H_4$, $H_5$, $H_6$, $H_7$ \cite{Hor89,HO31} and Lauricella's hypergeometric functions $F_A$, $F_B$, $F_C$, $F_D$~\cite{Lau93}.

In \cite{Beu13}, Beukers shows how to find a subgroup of the full monodromy group using Mellin--Barnes integral solutions of the associated $A$-hypergeometric system. This method only works under very restrictive conditions. These conditions are necessary to ensure the existence of a~basis of solutions in terms of Mellin--Barnes integrals. The monodromy groups found by Beukers' method are with respect to this Mellin--Barnes basis. In Sections~\ref{chap-ahyp}--\ref{chap-mon}, we will fix notation and introduce $A$-hypergeometric functions and Beukers' method.

The solution space of an $A$-hypergeometric system can be written in terms of integrals over generalised Pochhammer cycles \cite{Beu08} or over twisted cycles \cite{Kit94}. The unique invariant Hermitian form is the intersection form on these cycles. Actually computing this intersection form is computationally hard. Similarly it is in general computationally hard to explicitly compute an invariant hermitian form directly from a given monodromy group. In \cite{Beu89} Beukers and Heckman construct a unique invariant Hermitian form for the generalized hypergeometric system satisfied by ${}_{n+1}F_{n}$ using properties of its monodromy group. The goal of this paper is to give an explicit construction of the invariant Hermitian form over the monodromy group as constructed by Beukers' method \cite{Beu13}. The construction of this Hermitian form is given in Theorem \ref{MB-Hermitian} and its proof covers Sections~\ref{chap-herm} and~\ref{chap-res}. This invariant Hermitian form gives us insight into the structure of monodromy groups for $A$-hypergeometric functions with a Mellin--Barnes basis. The techniques provided in this paper may be used to extend these results to $A$-hypergeometric functions in general.

\section[The A-hypergeometric system]{The ${\boldsymbol A}$-hypergeometric system}\label{chap-ahyp}

Fix a positive integer $N$ and let $\bm{\gamma}\in\mathbb{R}^N$ be a row vector. Let $L\subset\mathbb{Z}^N$ be a lattice of
rank $d$ which satisfies the following conditions:
\begin{enumerate}\itemsep=0pt
\item $L$ is contained in the hyperplane $\big((l_1,l_2,\dots,l_N) \in \mathbb{Z}^N, \, \sum_{i=1}^Nl_i=0\big)$.
\item $L$ is saturated, i.e., $(L\otimes\mathbb{R})\cap\mathbb{Z}^N=L$.
\end{enumerate}
Now define
\begin{gather*}
\Phi_{\bm{\gamma}}^L:=\sum_{\bm{l}\in L}\prod_{j=1}^N\frac{z_j^{\gamma_j+l_j}}
{\Gamma(\gamma_j+l_j+1)}.
\end{gather*}
For the moment this is a formal series expansion. Throughout this paper we denote (row)-vectors with a bold face and its components have a normal font weight (e.g., $\bm{\gamma} = (\gamma_1,\dots,\gamma_N )$).
Notice that $\Phi_{\bm{\gamma}}^L
=\Phi_{\bm{\gamma}+\bm{l}}^L$ for any $\bm{l}\in L$.
Let $r = N-d$ and let $A$ be an $(r \times N)$-matrix with integer entries such that $L$ is
the integer kernel of $A$.
Let us define
$\bm{\alpha}=A\bm{\gamma}^\intercal$. Notice that $A(\bm{\gamma}+\bm{l})^\intercal
=\bm{\alpha}$ for any $\bm{l}\in L\otimes\mathbb{R}$.
We call this the parameter vector of the $A$-hypergeometric system we will define.
Because $L$ is contained in the hyperplane $\sum_{i=1}^Nl_i=0$, there is a linear form $h\colon \mathbb{R}^{r} \rightarrow \mathbb{R}$, where $h(\bm{a}) = 1$ for all column vectors $\bm{a}$ of $A$. A Gale dual of $A$, is an integer $d\times N$ matrix whose rows form a
$\mathbb{Z}$-basis of $L$, we denote this matrix by $B$.

It turns out that $\Phi_{\bm{\gamma}}^L$ satisfies a system of partial differential equations. First of all, let $\bm{m}=
(m_1,\dots,m_N)$ be an integer row vector such that $\bm{m}\cdot\bm{l}=0$ for all $\bm{l}\in L$.
Then, for any $\lambda \in \mathbb{C}^*$, one easily sees that
\begin{gather*}
\Phi_{\bm{\gamma}}^L(\lambda^{m_1}z_1,\dots,\lambda^{m_N}z_N)=
\lambda^{\bm{m}\cdot\bm{\gamma}}\Phi_{\bm{\gamma}}^L(z_1,\dots,z_N).
\end{gather*}
Take the derivative with respect to $\lambda$ and set $\lambda=1$. Then we see that
$\Phi_{\bm{\gamma}}^L$ is annihilated by the differential operator
\begin{gather*}
m_1z_1\partial_{z_1}+\cdots+m_Nz_N\partial{z_N}-\bm{m}\cdot\bm{\gamma}.
\end{gather*}
In particular, if we let $\bm{m}$ be the $i$-th row of $A=(A_{ij})$ we see
that $\Phi_{\bm{\gamma}}^L$ is annihilated by the Euler operator
\begin{gather*}
Z_i:=A_{i1}z_1\partial_{z_1}+\cdots+A_{iN}z_N\partial{z_N}-\alpha_i.
\end{gather*}
There is a second set of differential equations which arises from the
observation
\begin{gather*}
\partial_{z_1}^{\lambda_1}\cdots\partial_{z_N}^{\lambda_N}\Phi_{\bm{\gamma}}^L
=\Phi_{\bm{\gamma}-\bm{\lambda}}^L
\end{gather*}
for any $\bm{\lambda}=(\lambda_1,\dots,\lambda_N)\in\mathbb{Z}_{\ge0}^N$.
Let now $\bm{\lambda}\in L$ and write $\bm{\lambda}=\bm{\lambda}^+-\bm{\lambda}^-$,
where $\bm{\lambda}^\pm$ are integer vectors with non-negative entries. Then,
\begin{gather*}
\partial^{\bm{\lambda}^+}\Phi_{\bm{\gamma}}^L
=\Phi_{\bm{\gamma}-\bm{\lambda}^+}^L=\Phi_{\bm{\gamma}-\bm{\lambda}^-}^L
=\partial^{\bm{\lambda}^-}\Phi_{\bm{\gamma}}^L.
\end{gather*}
We use the notation $\partial^{\bm{\lambda}}=\partial_{z_1}^{\lambda_1}
\cdots\partial_{z_N}^{\lambda_N}$ and the second step follows from the invariance
of $\Phi_{\bm{\gamma}}^L$ when $\bm{\gamma}$ is shifted over vectors in $L$.
Thus we find that $\Phi_{\bm{\gamma}}^L$ is annihilated by the so-called box operators
\begin{gather*}
\Box^{\bm \lambda}:=\prod_{\lambda_i>0}\partial_{z_i}^{\lambda_i}-\prod_{\lambda_i<0}
\partial_{z_i}^{-\lambda_i}
\end{gather*}
for all $\bm{\lambda}\in L$.

{\samepage
The {\it $A$-hypergeometric system} $H_{A}(\bm{\alpha})$ is the system of differential equations generated by
\begin{enumerate}\itemsep=0pt
\item The Euler operators
\begin{gather*}
Z_j=A_{j1}\partial_{z_1}+\cdots+A_{jN}\partial_{z_N}-\alpha_j,\qquad j=1,\dots,N-d.
\end{gather*}
\item The box operators
\begin{gather*}
\Box^{\bm{\lambda}}=\partial^{\bm{\lambda}^+}-\partial^{\bm{\lambda}^-},\qquad
\bm{\lambda}\in L.
\end{gather*}
\end{enumerate}}

An $A$-hypergeometric function is a holomorphic
function in $z_1,\dots,z_N$ which satisfies the equations in the $A$-hypergeometric system. Either $A$ together with a parameter vector $\bm{\alpha}$ or~$B$
with $\bm{\gamma}$ (modulo translations in $L$) is enough to encode all the information about the $A$-hypergeometric system. The columns of $A$ are denoted by $\bm{a}_1,\dots,\bm{a}_N$ and the columns of $B$ are denoted by $\bm{b}_1,\dots,\bm{b}_N$.

\begin{Example}
Appell's $F_4$ is the hypergeometric function defined by
\begin{gather}
\HG{F_4}{a,\,\,b}{c,\,\,c'}{x,y} = \sum_{m=0}^{\infty}\sum_{n=0}^{\infty} \frac{(a)_{m+n}(b)_{m+n}}{(c)_m(c')_{n}m!n!}x^my^n.\label{eq-F4}
\end{gather}
Here $(a)_n$ denotes the Pochhammer symbol, defined as
\begin{gather*}
(a)_n = a(a+1)\cdots(a+n-1).
\end{gather*}
We can write Pochhammer symbols in terms of gamma functions as $(a)_n = \Gamma(a+n)/\Gamma(a)$. Recall Euler's reflection formula
\begin{gather*}
\Gamma(z)\Gamma(1-z) = \sin(\pi z).
\end{gather*}
Using this we can bring $\Gamma$-functions in the numerator in each summand of \eqref{eq-F4} down to the denominator. Then up to a constant factor we get
\begin{gather*}
\sum_{m=0}^{\infty}\sum_{n=0}^{\infty} \frac{x^my^n}{\Gamma(1-a-m-n)\Gamma(1-b-m-n)\Gamma(c+m)\Gamma(c'+n)\Gamma(m+1)\Gamma(n+1)}.
\end{gather*}
Substitute $x = \frac{z_3z_5}{z_1z_2}$ and $y = \frac{z_4z_6}{z_1z_2}$ and premultiply with $z_1^{-a}z_2^{-b}z_3^{c-1}z_4^{c'-1}$ to get
\begin{gather*}
\sum_{m=0}^{\infty}\sum_{n=0}^{\infty}\frac{z_1^{-a-m-n}} {\Gamma(1\!-a\!-m\!-n)} \frac{z_2^{-b-m-n}}{\Gamma(1\!-b\!-m\!-n)} \frac{z_3^{c-1+m}} {\Gamma(c+m)}\frac{z_4^{c'-1+n}}{\Gamma(c'+n)} \frac{z_5^m}{\Gamma(m+1)}  \frac{z_6^{n}}{\Gamma(n+1)}.
\end{gather*}
Let $L$ be the lattice generated by $(-1,-1,1,0,1,0)$ and $(-1,-1,0,1,0,1)$ and let $\bm{\gamma} = (-a,-b,\allowbreak c-1,c'-1,0,0)$ then this summation equals $\Phi_{\bm\gamma}^L$. In other words Appell's hypergeometric function
$\HG{F_4}{a,\,\,b}{c,\,\,c'}{x,y}$ can be obtained from the solutions of the $A$-hypergeometric system, where
\begin{gather*}
B = \begin{pmatrix}-1&-1&1&0&1&0\\-1&-1&0&1&0&1\end{pmatrix}\!,\qquad
\bm\gamma = (-a,-b,c-1,c'-1,0,0).
\end{gather*}
From this information we can also obtain an $A$-matrix and an $\bm{\alpha}$-vector pair. The lattice $L$ is generated by the rows of $B$ and $L$ is the integer kernel of $A$. So for $A$ and $\bm{\alpha}$ we can do linear algebra to obtain
\begin{gather*}
A = \begin{pmatrix}1 & 0 & 0 & 0 & 1 & 1 \\
0 & 1 & 0 & 0 & 1 & 1 \\
0 & 0 & 1 & 0 & -1 & 0 \\
0 & 0 & 0 & 1 & 0 & -1 \end{pmatrix}\!,\qquad
\bm{\alpha} = A\bm{\gamma}^{\intercal} = (-a, -b, c - 1, c' - 1)^{\intercal}.
\end{gather*}
\end{Example}
The $A$-hypergeometric systems we are interested in are those that are irreducible and only depend on $\bm{\alpha}$ modulo $\mathbb{Z}$. We can achieve this by assuming the system is non-resonant.

\begin{Definition}
An $A$-hypergeometric system $H_{A}(\bm{\alpha})$ is called {\it non-resonant} if the boundary of the cone $C(A) := \langle \bm{a}_1, \dots,\bm{a}_N \rangle_{\mathbb{R}_{\geq 0}}$ does not intersect the translated lattice $\bm{\alpha} + \mathbb{Z}^r$.
\end{Definition}

\begin{Theorem}[{\cite[Theorem 2.11]{GKZ90}}]
A non-resonant {$A$-hypergeometric} system $H_A(\bm\alpha)$ is irreducible.
\end{Theorem}

For reasons that will become clear in the next section we also want $\bm{\alpha}$ to be totally non-resonant.

\begin{Definition}\label{def:total-nonres}\sloppy
An $A$-hypergeometric system $H_{A}(\bm{\alpha})$ is called {\it totally non-resonant} if for each \mbox{$(r-1)$}-independent subset $\{\bm{a}_{j_1},\dots,\bm{a}_{j_{r-1}}\}$ of the set of columns of $A$ we have that $\langle \bm{a}_{j_1}, \dots,\bm{a}_{j_{r-1}} \rangle_{\mathbb{R}_{\geq 0}}$ does not intersect the lattice $\bm{\alpha} + \mathbb{Z}^r$.
\end{Definition}

We will always assume that $\bm{\alpha}$ is chosen totally non-resonant in the remainder of this paper unless otherwise stated.

\begin{Theorem}[{\cite[Corollary 5.20]{Ado94}}]\label{thm:A-rank}
Let $Q(A)$ be the convex hull of the points $\bm{a}_1,\dots,\bm{a}_N$. If the system $H_{A}(\bm{\alpha})$ is non-resonant then the holonomic rank of $H_{A}(\bm{\alpha})$ is equal to $\Vol(Q(A))$. Here the volume $\Vol$ is normalized such that the standard $(r-1)$-simplex has volume $1$.
\end{Theorem}
Let us denote by $D$ the holonomic rank of the $A$-hypergeometric system. Therefore when the system is non-resonant, then $D = \Vol(Q(A))$.

\section{Power series solutions}
\label{chap-power}
Recall the formal power series expansion
\begin{gather*}\Phi_{\bm{\gamma}}^L =
\sum_{\bm{l}\in L} \frac{\bm{z}^{\bm{\gamma} + \bm{l}}}{\bm{\Gamma}(\bm{\gamma} + \bm{l} + \bm{1}) }. 
\end{gather*}
Here and throughout this paper we use the convention that for any vector $\bm{v} = \left(v_1,\dots,v_N\right)$ the entity $\bm{\Gamma}(\bm{v})$ denotes $\prod_{i=1}^N \Gamma(v_i)$ and $\bm{z}^{\bm{v}} = \prod_{i=1}^N z_i^{v_i}$. And here $\bm{1}$ is the ones-vector. For a scalar~$c$ and vector $\bm{v}$ we let $c^{\bm{v}} = \left(c^{v_1},\dots,c^{v_N}\right)$.
We have seen in the previous section that $\Phi_{\bm{\gamma}}^L$ satisfies the $A$-hypergeometric system
$H_A(\bm\alpha)$ with $\bm\alpha=A\bm\gamma^\intercal$. Notice that these equations do
not change if we shift $\bm\gamma$ by a vector in $L\otimes\mathbb{R}$. Hence we get
in principle an infinite dimensional space of formal solutions. However, we shall only
be interested in those shifts of $\bm\gamma$ that yield Puiseux series solutions with
a domain of convergence. They belong to the $D$-dimensional solution space mentioned in
Theorem \ref{thm:A-rank}.

The question is now how to determine these shifts.
To answer this question we will use that $1/\Gamma(x)$ is $0$ if $x \in \mathbb{Z}_{\leq 0}$. Another observation is that if we let a basis for $L$ be $\bm{l}_1,\dots,\bm{l}_d$, then even though we have $N$ variables $z_1,\dots,z_N$, effectively we are only using
$x_1 = \bm{z}^{\bm{l}_1},\dots,x_d = \bm{z}^{\bm{l}_d}$. In~this way we can rewrite $\Phi_{\bm{\gamma}}^L$ as
\begin{gather*}
\Phi_{\bm{\gamma}} = \bm{z}^{\bm{\gamma}}\sum_{\bm{k} \in \mathbb{Z}^d} \frac{\bm{x}^{\bm{k}}}{\bm{\Gamma}(\bm{\gamma} + \bm{k}B + \bm{1}) },
\end{gather*}
where $B$ is the $(d\times N)$-matrix with $\bm{l}_i$ as its $i$-th row and $\bm{k}$ is considered as a row-vector.

To describe the shifts of $\bm\gamma$ we fix $\bm\gamma_0$ such that $\bm\alpha
=A\bm\gamma_0^{\intercal}$ and parametrize all shifts by $\bm\gamma_0+\bm\mu B$, where $\bm\mu
\in\mathbb{R}^d$ is considered as a row vector. Since $\Phi_{\bm\gamma+\bm{l}}=\Phi_{\bm\gamma}$
for all $\bm{l}\in L$, we can restrict $\bm\mu$ to the domain $[0,1)^d$. We can
now rewrite $\Phi_{\bm\gamma}=\bm{z}^{\bm\gamma_0}\Psi_{\bm\mu}$, where
\begin{gather*}
\Psi_{\bm\mu}=\sum_{\bm{k}\in\mathbb{Z}^d}\frac{\bm{x}^{\bm{k}+\bm\mu}}
{\bm{\Gamma}(\bm{\gamma_0} + (\bm{k}+\bm\mu)B + \bm{1}) }.
\end{gather*}

Fix $I\subset\{1,2,\dots,N\}$ with cardinality $d$ and $\bm{b}_i$, $i\in I$ linearly
independent. We call such a~set a {\it cotriangle}, the reason being that
the vectors $\bm{a}_i$, $i\in I^c$ span a simplex (triangle) in $\mathbb{Z}^r$.
Then choose $\bm\mu\in[0,1)^d$ such that $\bm\gamma_0+\bm\mu B$ has integer components at the indices $i\in I$. Let~$B_I$ be the submatrix of $B$ consisting of the
columns $\bm{b}_i$, $i\in I$ and let $\bm\gamma_{0I}$ be the sub-rowvector of~$\bm{\gamma}_0$ consisting of the indices in $I$. Then we need to solve $\bm\gamma_{0I}+\bm\mu B_I\in\mathbb{Z}^d$ in $\bm\mu\in
[0,1)^d$. This comes down to counting the number of shifted integral points in the fundamental
parallelogram spanned by the rows of $B_I$. Clearly the number of solutions is independent of this shift and there are precisely $\Delta_I$ solutions, where $\Delta_I=|\det(B_I)|$. Having found such $\bm\mu$ we note that the sum in
the definition of $\Psi_{\bm\mu}$ is restricted to the domain $\gamma_{0,i}+(\bm{k}+\bm\mu)
\bm{b}_i\ge0$, $i\in I$. This is because $1/\Gamma(x+1)=0$ if $x$ is a negative integer.

Choose a point $\bm\rho$ in the interior of the positive cone spanned by the $\bm{b}_i$, $i\in I$.
Then $\Psi_{\bm\mu}$ converges at the points $\bm{x}$ with $|x_i|=t^{\rho_i}$ for
sufficiently small $t>0$. See \cite{Beu11} for the necessary estimates.
We call $\bm\rho$ a convergence direction.

It is conceivable that besides $I$ there is another index at which $\bm\gamma_0+\bm\mu B$ has
an integer component. Since
\begin{gather*}
\bm\alpha=A\bm\gamma^\intercal=A\bm\gamma^\intercal+AB^\intercal\bm\mu^\intercal,
\end{gather*}
this means that $\bm\alpha$ can be written as a linear combination of the vectors $\bm{a}_i$
with fewer than $r=N-d$ non-integral coordinates. By our assumption of total non-resonance,
see Definition~\ref{def:total-nonres}, this situation cannot occur. We conclude that $I$
is uniquely determined by $\bm\mu$.

\begin{Definition}
We call $\bm\mu\in[0,1)^d$ a {\it solution point} and denote the corresponding set $I$ by $I(\bm\mu)$. Its corresponding parameter vector is denoted by $\bm{\gamma}^{\bm{\mu}} := \bm\gamma_0+\bm\mu B$.
\end{Definition}

Let us reverse the situation and start with a convergence direction $\bm\rho\in\mathbb{R}^d$
not in the hyperplane spanned by any $d-1$ vectors $\bm{b}_i$. The set of cotriangles
$I$ such that $\bm\rho$ is contained in the positive cone generated by $\bm{b}_i$, $i\in I$
is denoted by $\mathcal{I}_{\bm\rho}$. Each cotriangle $I$ contributes $\Delta_I$ solution
points $\bm\mu$ and so we find $\sum_{I\in\mathcal{I}_{\bm\rho}}\Delta_I$ Laurent series
solutions that converge around $\bm\rho$. We call $\mathcal{I}_{\bm\rho}$ a
{\it cotriangulation} of $B$.

From \cite[Section 5.4]{DRS10} it follows that cotriangulations of $B$ are
in one-to-one correspondence with triangulations of $A$. The correspondence is given
by associating a cotriangle $I$ with a~triangle (simplex) spanned by $\bm{a}_i$, $i\in I^c$.
Furthermore, it follows from \cite[Lemma 14.2]{Beu11} that $\Delta_I=|\det(\bm{a}_i)_{i\in I^c}|$.
Hence $\sum_{I\in\mathcal{I}_{\bm\rho}}\Delta_I$ equals ${\rm Vol}(Q(A))$, which is
precisely the rank of our hypergeometric system $H_A(\bm\alpha)$. Thus the Laurent
series $\bm{z}^{\bm\gamma_0}\Psi_{\bm\mu}$ with $I(\bm\mu)\in\mathcal{I}_{\bm\rho}$
forms a basis of solutions with a common domain of convergence.

\begin{Definition}
\label{def-chamber-wall-fan}
A {\it chamber} is a fully dimensional cone constructed as an intersection of the form $\mathcal{C}_{\bm\rho} := \bigcap_{I \in \mathcal{I}_{\bm\rho}} C_I$, where $C_I$ is the cone generated by the $\bm{b}_i$ for $i \in I$. It has the property that for each convergence direction $\bm\rho'$ we pick in the interior of the chamber that $\mathcal{C}_{\bm\rho'} = \mathcal{C}_{\bm\rho}$. In~this way cotriangulations and chambers are in one-to-one correspondence. A {\it wall} is any face of a chamber that is of codimension $1$. The polyhedral complex generated by the chambers $\mathcal{C}_{\bm\rho}$ and all of their faces is called the {\it secondary fan} $\Sigma_B$.
\end{Definition}

\section{Mellin--Barnes integrals}

Let notation be as above and choose a vector $\bm{\sigma} \in \mathbb{R}^d$. For any vector $\bm{s} = (s_1,\dots,s_d)$ denote ${\rm d} \bm{s}= {\rm d} s_1 \wedge {\rm d} s_2 \wedge \dots \wedge {\rm d} s_d$. Then consider the integral
\begin{gather*}
M(\bm{z}) = M(z_1,\dots,z_N) := \int_{\bm{\sigma} + i\mathbb{R}^d} \bm{\Gamma}(-\bm{\gamma}_0-\bm{s}B)\bm{z}^{\bm{\gamma}_0 + \bm{s}B} {\rm d} \bm{s}.
\end{gather*}

This is a so-called Mellin--Barnes integral. When there is a basis of solutions
for an $A$-hypergeometric system in terms of Mellin--Barnes integrals, then this will help us to find the monodromy group for this $A$-hypergeometric system. A quick summary about Mellin--Barnes integrals is given here, for a more thorough introduction see \cite{Beu13}.

\begin{Theorem}[{\cite[Theorem 3.1]{Beu13}}]\label{thm-MB-gamma}
Suppose that $\gamma_{0,i} < - \bm{b}_i \cdot \bm{\sigma}$ for
$i = 1,\dots,N$ and that $M(\bm{z})$ converges. Then $M(\bm{z})$ satisfies the differential system $H_{A}(\bm{\alpha})$.
\end{Theorem}
Now not all systems admit a choice for $\bm{\gamma}_0$, where $\gamma_{0,i} < - \bm{b}_i \cdot \bm{\sigma}$. Using contiguity relations we can change $\bm{\alpha}$ without affecting the monodromy and we still have a freedom in $\bm{\sigma}$.
In \cite{Beu13} it is shown that we can choose $\bm{\sigma}$ and $\bm{\alpha}$ such that $\bm{\gamma}_0$ satisfies the conditions of Theorem \ref{thm-MB-gamma} without affecting monodromy.

For convergence of Mellin--Barnes solutions we will define the open zonotope
\begin{gather*}
Z_B^{\circ} = \Bigg\{\sum_{i=1}^N \nu_i \bm{b}_i\,\bigg|\,\, 0 < \nu_i < 1\bigg\}.
\end{gather*}
Note that our definition of a zonotope is scaled with a factor two compared to its definition in~{\rm \cite{Beu13}}.

\begin{Remark}
\label{rem-zonotope-sym}
The zonotope $Z_B^{\circ}$ is convex by definition and it is centrally symmetric. Namely if $v \in Z_B^{\circ}$ then $-v \in Z_B^{\circ}$. Indeed, there exists numbers $\nu_i \in (0,1)$ such that $v = \sum_{i=1}^N \nu_i\bm{b}_i$. As~the sum of the column vectors of $B$ is equal to $0$ this means \begin{gather*}
-v = \sum_{i=1}^{N} \bm{b}_i - v = \sum_{i=1}^N (1 - \nu_i)\bm{b}_i \in Z_B^{\circ}.
\end{gather*}
\end{Remark}
We introduce the variables $\bm{x}=\bm{z}^B$ and rewrite $M(\bm{z})$ as
$\bm{z}^{\bm\gamma_0}\tilde{M}(\bm{x})$, where
\begin{gather*}
\tilde{M}(\bm{x})=\int_{\bm{\sigma} + i\mathbb{R}^d} \bm{\Gamma}(-\bm{\gamma_0}-\bm{s}B)
\bm{x}^{\bm{s}} {\rm d} \bm{s}.
\end{gather*}
\begin{Theorem}[{\cite[Corollary 4.2]{Beu13}}] 
Let $\bm{\tau} = \frac{1}{2\pi}\Arg(\bm{x})$ be a component-wise choice of argument of the vector $\bm{x}$. Then $\tilde{M}(\bm{x})$ converges absolutely if
$\bm{\tau} \in \frac{1}{2}Z_B^{\circ}$.
\end{Theorem}

Lastly we quickly state how linearly independent solutions can be found, and thus how we can find a basis of solutions using Mellin--Barnes integrals. The following theorem tells us that choosing different $\bm\tau\in\frac12 Z_B^\circ$ we can obtain independent Mellin--Barnes solutions.

\begin{Theorem}[{\cite[Proposition 4.6]{Beu13}}] 
Let $H_A(\bm{\alpha})$ be a non-resonant $A$-hypergeometric system. Let
$\bm{\tau}_1,\dots,\bm{\tau}_q \in \frac{1}{2}Z_B^{\circ}$ be points whose coordinates
differ by integers. Fix a point $\bm{x}^0\in(\mathbb{C}^\times)^d$
and choose for each $\bm\tau_i$ the
Mellin--Barnes integral $\tilde{M}_i(\bm{x})$ with this argument choice for $\bm{x}^0$.
Then $\tilde{M}_1,\dots,\tilde{M}_q$ are linearly independent in a neighbourhood of $\bm{x}^0$.
\end{Theorem}
In particular this implies that if $q=D$, then we have a basis of solutions of
$H_A(\bm\alpha)$ given by Mellin--Barnes integrals.

\section{Monodromy}
\label{chap-mon}
In this section we consider an $A$-hypergeometric system with solution space $V$. Any non-zero solution in $V$ remains a non-zero solution if we analytically continue it around some loop.
This means that analytic continuation along some loop $c$ induces a linear map
$\phi_{c}\colon V \rightarrow V$.
All the possible elements $\phi_{c}$ give the monodromy group.
Seeing the elements $\phi_{c}$ as matrices, then the monodromy group will depend on a choice of basis. In our case this basis will be given by Mellin--Barnes solutions $\tilde{M}_1,\dots,\tilde{M}_D$.

Let $\bm{n} \in \mathbb{Z}^d$ be a column vector and let $c(\bm{n})$ be the loop
\begin{gather*}
\big\{\big({\rm e}^{2\pi {\rm i} n_1 t}x_1,\dots,{\rm e}^{2\pi {\rm i} n_d t}x_d\big)\mid t \in [0,1]\big\}.
\end{gather*}
Analytic
continuation of the Laurent series solution $\Psi_{\bm{\mu}}(\bm{x})$ along $c(\bm{n})$
gives ${\rm e}^{2\pi {\rm i} \bm{n} \cdot \bm\mu}\Psi_{\bm\mu}(\bm{x})$. This means that, given a convergence direction $\bm\rho$, and its corresponding basis of local Laurent series solutions
$\Psi_{\bm\mu_1},\dots,\Psi_{\bm\mu_D}$, the monodromy elements $\phi_{c(\bm{n})}$ can be written in matrix form as
\begin{gather*}
\chi_{\bm\rho,\bm{n}} := \begin{pmatrix}{\rm e}^{2\pi {\rm i} \bm\mu_1 \bm{n}} & 0 & \cdots & 0 \\ 0 & {\rm e}^{2\pi {\rm i} \bm\mu_2 \bm{n}} & \cdots & 0 \\ \vdots & \vdots & \ddots & \vdots \\ 0 & 0 & \cdots & {\rm e}^{2\pi {\rm i} \bm\mu_D \bm{n}}\end{pmatrix}\!.
\end{gather*}

This gives a commutative subgroup of the monodromy group which is generated by the elements $\chi_{\bm\rho,j} := \chi_{\bm\rho,\bm{e}_j}$, $j=1,\dots,d$.

Now suppose that $H_{A}(\bm{\alpha})$ has a Mellin--Barnes basis of solutions and therefore
there exists a~set $\bm{\tau}_1,\dots,\bm{\tau}_D$ such that
$\bm{\tau}_i \in \frac{1}{2}Z_B^{\circ}$ are distinct and differ by integers.
Denote the Mellin--Barnes integral corresponding to the argument choice
$2\pi\bm\tau_j$ by $\tilde{M}_j$.

Consider the Mellin--Barnes basis near a point $\bm{x}^0$. Analytic continuation
of $\tilde{M}_1$ along the path $c(\bm\tau_j-\bm\tau_1)$ changes $\tilde{M}_1$ into $\tilde{M}_j$. Note that this
is independent of the choice of $\bm{x}^0$.
If we write a local series expansion
$\tilde{M}_1 = \sum_{k=1}^{D} \lambda_k \Psi_{\bm\mu_k}$ for some
convergence direction $\bm\rho$, then analytic continuation along
$c(\bm\tau_j-\bm\tau_1)$ will result in
$\tilde{M}_j = \sum_{k=1}^{D} \lambda_k {\rm e}^{2\pi {\rm i}(\bm\tau_j-\bm\tau_1)\cdot \bm\mu_k}
\Psi_{\bm\mu_k}$.
If one of these $\lambda_k$'s is zero, we see that $\tilde{M}_1,\dots,\tilde{M}_D$ spans a space of dimension strictly less than $D$, which is in contradiction with $\tilde{M}_1,\dots,\tilde{M}_D$ being linearly independent. Hence it must be that the $\lambda_k$'s are all non-zero. We can then normalize the $\Psi_{\bm\mu_k}$ such that the $\lambda_k$'s are $1$ and obtain a transition matrix between Mellin--Barnes solutions to local power series solutions:
\begin{gather}X_{\bm\rho} = \begin{pmatrix}
 1 & 1 & \cdots & 1 \\
 {\rm e}^{2\pi {\rm i} \bm\mu_1(\bm\tau_{2} - \bm\tau_{1})} & {\rm e}^{2\pi {\rm i} \bm\mu_2(\bm\tau_{2} - \bm\tau_{1})} & \cdots & {\rm e}^{2\pi {\rm i} \bm\mu_D(\bm\tau_{2} - \bm\tau_{1})} \\
 {\rm e}^{2\pi {\rm i} \bm\mu_1(\bm\tau_{3} - \bm\tau_{1})} & {\rm e}^{2\pi {\rm i} \bm\mu_2(\bm\tau_{3} - \bm\tau_{1}) } & \cdots & {\rm e}^{2\pi {\rm i} \bm\mu_D(\bm\tau_{3} - \bm\tau_{1})} \\
 \vdots & \vdots & \ddots & \vdots \\
 {\rm e}^{2\pi {\rm i} \bm\mu_1(\bm\tau_{D} - \bm\tau_{1})} & {\rm e}^{2\pi {\rm i} \bm\mu_2(\bm\tau_{D} - \bm\tau_{1})} & \cdots & {\rm e}^{2\pi {\rm i} \bm\mu_D(\bm\tau_{D} - \bm\tau_{1})}
 \end{pmatrix}\!,\label{eq-Xrho}\end{gather}
which satisfies
\begin{gather*}
\begin{pmatrix}
\tilde{M}_1\\ \tilde{M}_2\\ \tilde{M}_3\\ \vdots \\ \tilde{M}_D\end{pmatrix} = X_{\bm\rho} \begin{pmatrix}\Psi_{\bm\mu_1} \\ \Psi_{\bm\mu_2} \\ \Psi_{\bm\mu_3} \\ \vdots \\ \Psi_{\bm\mu_D}
\end{pmatrix}\!.
\end{gather*}

This means that the monodromy subgroup generated by $\chi_{\bm\rho,j}$ with respect to a basis of local series expansions, can be transformed through $X_{\bm\rho}$ into a monodromy subgroup with respect to a basis of Mellin--Barnes solutions.

The matrices that generate this monodomy subgroup with respect to a basis of Mellin--Barnes solutions are defined as
\begin{gather*}
M_{\rho,j} = X_{\bm\rho}\chi_{\bm\rho,j}X_{\bm\rho}^{-1}.
\end{gather*}

By changing the convergence direction $\bm\rho$ we will therefore obtain multiple subgroups, which together will generate a larger subgroup of the monodromy group $\mathcal{M}$. Since it is unclear whether this generates the whole monodromy group, we will define a subgroup of the monodromy.

\begin{Definition}
The Mellin--Barnes group $\mathcal{M}_{MB}$ is the group generated by the matrices $M_{\bm\rho,j}$ for all $j=1,\dots,d$ and convergence directions $\bm\rho$.
\end{Definition}

In \cite{Beu13} it is shown that the Mellin--Barnes group and the monodromy group are equal for the systems $F_2$ and ${}_{n+1}F_n$. However, in general it is not easy to check whether the monodromy group and Mellin-barnes group are equal. The reason for this is that most computations of known monodromy groups use a basis that is hard to translate in terms of Mellin--Barnes integrals.
\begin{Remark}
The Mellin--Barnes group corresponds to the power series $\Psi_{\bm\mu}$, though we started out with the power series $\Phi_{\bm\gamma}$. These power series differ by a monomial factor. Hence their corresponding monodromy groups are the same upto multiplication by scalars.
\end{Remark}

\section{The Hermitian form}\label{chap-herm}
In this section we adopt the notations from the sections above. In particular $X_{\bm\rho}$ are the transition matrices given in \eqref{eq-Xrho}. Our goal is to prove the following theorem.

\begin{Theorem}\label{MB-Hermitian}
Let $H_A(\bm{\alpha})$ be a totally non-resonant $A$-hypergeometric system admitting a~Mel\-lin--Barnes basis of solutions. Then there exists a non-trivial Hermitian form $H$ which is inva\-riant under the group $\mathcal{M}_{MB}$. Furthermore given any convergence direction $\bm\rho$, this Hermitian form can be given explicitly as
\begin{gather*}
H = \big(\bar{X}_{\bm\rho}^{\intercal}\big)^{-1} \Delta_{\bm\rho} X_{\bm\rho}^{-1},
\end{gather*}
where $\Delta_{\bm\rho}$ is the diagonal matrix
\begin{gather}
\operatorname{diag}\bigg(\bigg\{\Delta_{I_k} \prod_{l \in I_k} (-1)^{\gamma_l^{\bm\mu_k}}\prod_{i \not \in I_k} \sin\big(\pi \gamma_i^{\bm\mu_k}\big)\bigg\}_{k=1,\dots,D}\bigg)
\label{eq-deltarho}
\end{gather}
and where $\bm\mu_k$ runs over all solution points with $I_k := I(\bm{\mu}_k) \in \mathcal{I}_{\bm\rho}$.
\end{Theorem}

\begin{Note}
Due to lack of space for certain formulas and equations, we sometimes use a different notation for matrices. In our case for an $M\times N$ matrix, where $M$ and $N$ are known we use the notation \begin{gather*}
\left\{a_{rc}\right\}_{r,c} := \begin{pmatrix}a_{11} & a_{12} & \cdots & a_{1N} \\ a_{21} & a_{22} & \cdots & a_{2N} \\ \vdots & \vdots & \ddots & \vdots \\ a_{M1} & a_{M2} & \cdots & a_{MN} \end{pmatrix}\!.
\end{gather*}
For diagonal matrices of fixed dimension $N$ we may use the notation
\begin{gather*}
\left\{a_{r}\right\}_{rr} := \begin{pmatrix}a_{1} & 0 & \cdots & 0 \\ 0 & a_{2} & \cdots & 0 \\ \vdots & \vdots & \ddots & \vdots \\ 0 & 0 & \cdots & a_{N} \end{pmatrix}\!.
\end{gather*}
\end{Note}

\begin{proof}
Fix a convergence direction $\bm\rho$ and consider $H_{\bm\rho} = \big(\bar{X}_{\bm\rho}^{\intercal}\big)^{-1} \Delta_{\bm\rho} X_{\bm\rho}^{-1}$, where $X_{\bm{\rho}}$ is the transitition matrix given in \eqref{eq-Xrho} and $\Delta_{\bm\rho}$ is given in \eqref{eq-deltarho}. We show that $H_{\bm\rho}$ is an invariant Hermitian matrix for the monodromy matrices $M_{\bm\rho,j} = X_{\bm\rho} \chi_{\bm\rho,j} X_{\bm\rho}^{-1}$ defined in Section \ref{chap-mon}. It is clear that $H_{\bm\rho}$ is an Hermitian matrix. Next we show that $H_{\bm\rho}$ is invariant under Hermitian conjugation with $M_{\bm\rho,j}$, thus we want
\begin{gather*}
\big(\overline{X_{\bm\rho} \chi_{\bm\rho,j} X_{\bm\rho}^{-1}}\big)^{\intercal}\big(\bar{X}_{\bm\rho}^{\intercal}\big)^{-1} \Delta_{\bm\rho} X_{\bm\rho}^{-1} X_{\bm\rho} \chi_{\bm\rho,j} X_{\bm\rho}^{-1} = \big(\bar{X}_{\bm\rho}^{\intercal}\big)^{-1} \Delta_{\bm\rho} X_{\bm\rho}^{-1}.
\end{gather*}
We can simplify this into
\begin{gather*}
\overline{\chi_{\bm\rho,j}}^{\intercal}\Delta_{\bm\rho} \chi_{\bm\rho,j} = \Delta_{\bm\rho}.
\end{gather*}
As all of these matrices are diagonal, and $\overline{\chi_{\bm\rho,j}}^{\intercal}$,
$\chi_{\bm\rho,j}$ are each others inverse we see that this is indeed an equality.

The remainder of the proof consists of showing that $H_{\bm\rho}$ is independent of the
choice of $\bm\rho$. The resulting matrix $H$ is then an invariant Hermitian form for all
local monodromy matrices~$M_{\bm\rho,j}$.

As explained in Section \ref{chap-power} we associate to each convergence direction a set of solution points $\bm{\mu}_1,\dots,\bm{\mu}_D$ and cotriangles $I_k := I(\bm{\mu}_k).$ To prove the independence of $H_{\bm\rho}$,
we calculate $H^{-1}_{\bm\rho}$, where we denote $\tilde{\bm{\tau}}_l = \bm{\tau}_l - \bm{\tau}_1$,
\begin{align*}
H_{\bm\rho}^{-1}&= X_{\bm\rho} \Delta_{\bm\rho}^{-1} \bar{X}_{\bm\rho}^{\intercal}\\
&= \big\{{\rm e}^{2\pi {\rm i}\bm{\mu}_c \tilde{\bm{\tau}}_r}\big\}_{r,c}
\Bigg\{\frac{1}{\Delta_{I_r}} \prod_{l \in I_r} (-1)^{\gamma_l^{\bm\mu_r}}\prod_{i \not \in I_r} \frac{1}{\sin(\pi \gamma_i^{\bm\mu_r})}\Bigg\}_{rr}
\big\{{\rm e}^{-2\pi {\rm i}\bm{\mu}_r \tilde{\bm{\tau}}_c}\big\}_{r,c}
\\
&= \Bigg\{\frac{{\rm e}^{2\pi {\rm i}\bm{\mu}_c \tilde{\bm{\tau}}_r}}{\Delta_{I_c}} \prod_{l \in I_c} (-1)^{\gamma_l^{\bm\mu_c}}\prod_{i \not \in I_c} \frac{1}{\sin(\pi \gamma_i^{\bm\mu_c})}\Bigg\}_{r,c}
\big\{{\rm e}^{-2\pi {\rm i}\bm{\mu}_r \tilde{\bm{\tau}}_c}\big\}_{r,c}
\\
&= \Bigg\{\sum_{k=1}^D \frac{{\rm e}^{2\pi {\rm i}\bm{\mu}_k (\tilde{\bm{\tau}}_r - \tilde{\bm{\tau}}_c)}\prod_{l \in I_k} (-1)^{\gamma_l^{\bm\mu_k}}}{\Delta_{I_k}\prod_{i \not \in I_k} \sin(\pi \gamma_i^{\bm\mu_k})} \Bigg\}_{r,c}\\
&=(2{\rm i})^r\Bigg\{\sum_{k=1}^D \frac{{\rm e}^{2\pi {\rm i}\bm{\mu}_k (\bm{\tau}_r - \bm{\tau}_c)}}{\Delta_{I_k}}\prod_{l \in I_k} {\rm e}^{\pi {\rm i} \gamma_l^{\bm\mu_k}}\prod_{l \not \in I_k} \frac{1}{{\rm e}^{\pi {\rm i} \gamma_l^{\bm\mu_k}} - {\rm e}^{-\pi {\rm i} \gamma_l^{\bm\mu_k}}} \Bigg\}_{r,c} \\
&= (2{\rm i})^r \Bigg\{\sum_{k=1}^D \frac{{\rm e}^{2\pi {\rm i}\bm{\mu}_k (\bm{\tau}_r - \bm{\tau}_c)}}{\Delta_{I_k}}\prod_{l = 1}^N {\rm e}^{\pi {\rm i} \gamma_l^{\bm\mu_k}}\prod_{l \not \in I_k} \frac{1}{{\rm e}^{2\pi {\rm i} \gamma_l^{\bm\mu_k}} - 1} \Bigg\}_{r,c} \\
&= (2{\rm i})^r \prod_{j=1}^N {\rm e}^{\pi {\rm i} \gamma_{0j}}\Bigg\{\sum_{k=1}^D \frac{{\rm e}^{2\pi {\rm i}\bm{\mu}_k (\bm{\tau}_r - \bm{\tau}_c)}}{\Delta_{I_k}}\prod_{l \not \in I_k} \frac{1}{{\rm e}^{2\pi {\rm i} \gamma_l^{\bm\mu_k}} - 1} \Bigg\}_{r,c}.
\end{align*}
Here the last step follows from the fact that the sum of the column vectors of $B$ is $0$ and $\bm\gamma_l^{\bm \mu_k} = \bm\gamma_0 + \bm{\mu}_kB$.
Each component of the inner matrix will be linked to a sum of certain residues, which can be seen from Lemma \ref{eval-res-complete} below. Using this and using
$\bm{\tau}_{r} - \bm{\tau}_{c} \in Z_B^{\circ}$ (see Remark~\ref{rem-zonotope-sym}) it follows from Corollary \ref{cor-gen-tri} below that $H_{\bm\rho}$ is independent of the choice
of $\bm\rho$.
\end{proof}

\section{Residues}
\label{chap-res}
Define the following differential form
\begin{gather*}
\omega := \omega(\bm{\tau},\bm{z}) = \frac{\bm{z}^{\bm\tau}}{(x_1\bm{z}^{\bm{b}_1} - 1)(x_2\bm{z}^{\bm{b}_2} - 1) \cdots (x_N\bm{z}^{\bm{b}_N} - 1)} \frac{{\rm d}\bm{z}}{\bm{z}},
\end{gather*}
where $x_j={\rm e}^{2\pi {\rm i}\gamma_{0,j}}$.
Here $\frac{{\rm d}\bm{z}}{\bm{z}}$ is short for $\frac{{\rm d}z_1}{z_1} \wedge \dots \wedge \frac{{\rm d}z_d}{z_d}$. And $\bm{z}^{\bm{b}}$ stands for $z_1^{b_1}\cdots z_d^{b_d}$.
Certain residues of this form are special cases of so-called {\it binomial residues} \cite{Cat00}.
For any solution point~$\bm{\mu}$, define the vector
\begin{gather*}
\bm{\zeta}^{\bm{\mu}} := {\rm e}^{2\pi {\rm i} \bm{\mu}},
\end{gather*}
where we use the notation ${\rm e}^{2\pi {\rm i}\bm{v}}=\big({\rm e}^{2\pi {\rm i}v_1},\dots,{\rm e}^{2\pi {\rm i}v_d}\big)$.

Notice that
\begin{gather*}
x_i (\bm{\zeta}^{\bm{\mu}})^{\bm{b}_i}={\rm e}^{2\pi {\rm i}\gamma_{0,i}}
{\rm e}^{2\pi {\rm i}\bm{\mu}\bm{b}_i}={\rm e}^{2\pi {\rm i}\gamma^{\bm{\mu}}_i}=1,
\end{gather*}
for all $i\in I(\bm{\mu})$ because $\gamma^{\bm\mu}_i\in\mathbb{Z}$ for all $i\in I(\bm\mu)$.
We thus see that $\bm{\zeta}^{\bm\mu}$ is a solution to the system of equations
$x_i\bm{z}^{\bm{b}_i}-1=0$, $i\in I(\bm\mu)$ in $\bm{z}$.

Let $f_i = x_i\bm{z}^{\bm{b}_i} - 1$ for $i=1,\dots,N$. Following \cite[p.~650]{GH78} we may define the residue
\begin{gather*}
\Res_{\bm{z} = \bm{\zeta}^{\bm\mu}} \omega = \pm \frac{(\bm{\zeta}^{\bm\mu})^{\bm{\tau}}}{J_{I}(\bm{\zeta}^{\bm\mu})\prod_{j \in I^c} f_j(\bm{\zeta}^{\bm\mu}) },
\end{gather*}
for $I = I(\bm{\mu}) = (i_1,\dots,i_d)$, where we choose the sign $\pm$ to be $\sign(\det(B_I))$
and where $J_I$ is the toric Jacobian given by
\begin{gather*}
J_I = \det\bigg(\bigg\{z_r \frac{\partial f_{i_c}}{\partial z_r}\bigg\}_{r,c}\bigg).
\end{gather*}
The choice of sign is not usually part of the definition of a residue. It is used in this context to make the value of the residue orientation independent.
Due to the simplicity of the functions $f_i$ we can easily show that
\begin{gather*}
J_I = \det(B_I) \prod_{j \in I}x_j\bm{z}^{\bm{b}_{j}}.
\end{gather*}
By definition of $\bm{\zeta}^{\bm\mu}$ we get $\prod_{j \in I}x_j(\bm{\zeta}^{\bm\mu})^{\bm{b}_{j}} = 1$, so as a consequence we get
\begin{gather}
\Res_{\bm{z} = \bm{\zeta}^{\bm\mu}} \omega =
\frac{(\bm{\zeta}^{\bm\mu})^{\bm{\tau}}}{\Delta_I
\prod_{j \in I^c} f_j(\bm{\zeta}^{\bm\mu})}.
\label{eq-res-omg}
\end{gather}
\noindent Lemma \ref{eval-res-complete} is now a direct consequence of \eqref{eq-res-omg}.

\begin{Lemma}\label{eval-res-complete}
Let $\bm\mu$ be a solution point then we have
\begin{gather*}
\Res_{\bm{z} = \bm{\zeta}^{\bm\mu}} \omega(\bm\tau,\bm{z}) = \frac{{\rm e}^{2 \pi {\rm i} \bm{\mu} \bm\tau}}{\Delta_I\prod_{j \in I^c}\big({\rm e}^{2 \pi {\rm i} \gamma^{\bm\mu}_j} - 1\big)},
\end{gather*}
where $I = I(\bm\mu)$.
\end{Lemma}

Using these residues we can now write a typical entry of the matrix
$H_{\bm\rho}^{-1}$ in the proof of Theorem \ref{MB-Hermitian} as
\begin{gather*}
\sum_{\bm\mu\colon I(\bm\mu)\in\mathcal{I}_{\bm\rho}}\Res_{\bm{z}=\bm{\zeta}^{\bm\mu}}
\omega(\bm{\tau}_r-\bm{\tau}_c,\bm{z}).
\end{gather*}
It would be tempting to prove that such an entry is independent of $\bm\rho$,
and hence the correspon\-ding cotriangulation $\mathcal{I}$, by
using general properties of multidimensional residues. Unfortunately we have
been unable to do so.
Instead we shall follow a local approach where we show equality of these sums
for neighbouring cotriangulations. In doing so we shall make use of residue
calculus for one variable rational functions.

Recall Definition \ref{def-chamber-wall-fan}.
\begin{Definition}
For any wall $W$ of the chamber $C_{\mathcal{I}}$ we denote by $\mathcal{I}_W$ all the cotriangles $I \in \mathcal{I}$ whose cones $C_I$ have $W$ as a (sub)-face.
\end{Definition}

\begin{Definition}
Two cotriangulations $\mathcal{I}$ and $\mathcal{J}$ are called adjacent if their corresponding chambers share the same wall. We call this wall the common wall between $\mathcal{I}$ and $\mathcal{J}$.
\end{Definition}

Given adjacent cotriangulations $\mathcal{I}$ and $\mathcal{J}$ with common wall $W$ then a cotriangle $I \in \mathcal{I}_W$ is characterized by having $d-1$ indices $i_1,\dots,i_{d-1}$ for which the cone generated by $\bm{b}_{i_1},\dots,\bm{b}_{i_{d-1}}$ contains $W$. The remaining index of $I$ corresponds to a $\bm{b}_{i_d}$ being on either side of $W$. Conversely, if we are given indices $i_1,\dots,i_{d-1}$
for which the corresponding cone generated by $\bm{b}_{i_1},\dots,\bm{b}_{i_{d-1}}$ contains
$W$ and if we are given an index $i_d$ for which the $\bm{b}$-vector is not on the hyperplane
$\operatorname{Hyp}(W)$, then $I = (i_1,\dots,i_{d})$ is either in $\mathcal{I}_W$ or
$\mathcal{J}_W$, depending on which side of the wall $\bm{b}_{i_d}$ lies.

\begin{Proposition}\label{lem-edge-res}
Let $\mathcal{I}$ and $\mathcal{J}$ be two adjacent cotriangulations with common wall $W$ and suppose $\bm\tau \in Z_B^\circ$ then
\begin{gather*}
\sum_{\bm\mu\colon I(\bm{\mu})\in\mathcal{I}_W}\Res_{\bm{z}=\bm\zeta^{\bm\mu}}\omega(\bm{\tau},\bm{z})
=\sum_{\bm\nu\colon I(\bm{\nu})\in\mathcal{J}_W}\Res_{\bm{z}=\bm\zeta^{\bm\nu}}\omega(\bm{\tau},\bm{z}).
\end{gather*}
\end{Proposition}

\begin{proof}
Choose any $i_1,\dots,i_{d-1}$ such that $\bm{b}_{i_1},\dots,\bm{b}_{i_{d-1}}$ are linearly independent and the cone spanned by them contains $W$. It suffices to prove
our lemma in case the sums run over all $I\in\mathcal{I}_W$, $J\in\mathcal{J}_W$
which contain $i_1,\dots,i_{d-1}$. The full lemma then follows after summation over
all sets $i_1,\dots,i_{d-1}$ such that the cone spanned by $\bm{b}_{i_1},\dots,
\bm{b}_{i_{d-1}}$ contains $W$.

Choose coordinates in $\mathbb{Z}^d$ such that the $d$-th coordinates of
$\bm{b}_{i_1},\dots,\bm{b}_{i_{d-1}}$ are zero. In gene\-ral we denote the $d$-th
coordinate of $\bm{b}_i$ by $\beta_i$. Hence $\beta_i=0$ for $i=i_1,\dots,i_{d-1}$.
Write $\bm{z}^{\bm{b}_i} = Q_i(z_1,\dots,z_{d-1})z_d^{\beta_i}$, where $Q_i$ is a monomial in $z_1,\dots,z_{d-1}$. Similarly we write $\bm{z}^{\bm\tau}=Q_0(z_1,\dots,z_{d-1})
z_d^{\tau_d}$.
Let $\delta$ be the determinant of $\big(\bm{b}_{i_1},\dots,\bm{b}_{i_{d-1}}\big)$
where we remove the last row, which is zero.
Then by construction we have that for any $i$ the following holds
\begin{gather*}
\det\left(\bm{b}_{i_1},\dots,\bm{b}_{i_{d-1}},\bm{b}_i\right) = \beta_i \delta.
\end{gather*}
The sign of $\beta_i$ determines on which side of $W$ the vector $\bm{b}_i$ lies.
Choose an index $i_d$ with $\beta_{i_d}\ne0$ and let $\bm{\zeta}$ be a point such that
$x_j\bm{\zeta}^{\bm{b}_j} = 1$ for $j \in I:=\{i_1,\dots,i_d\}$. Then \eqref{eq-res-omg} tells us that
\begin{gather}\label{res-at-point}
\Res_{\bm{z}=\bm\zeta}\omega(\bm\tau,\bm{z})
=\frac{Q_0(\zeta_1,\dots,\zeta_{d-1})\zeta_d^{\tau_d}}
{\Delta_I\prod_{j \not \in \{i_1,\dots,i_d\}}
\big(x_jQ_j(\zeta_1,\dots,\zeta_{d-1})\zeta_d^{\beta_j} - 1\big)}.
\end{gather}
We like to write this as a one variable residue. The variable will be called $w$. Consider
\begin{gather*}
\Omega(w)=\frac{Q_0(\zeta_1,\dots,\zeta_{d-1})w^{\tau_d}}
{\prod_{j\not\in\{i_1,\dots,i_{d-1}\}}
\big(x_jQ_j(\zeta_1,\dots,\zeta_{d-1})w^{\beta_j}-1\big)}\frac{{\rm d}w}{w}.
\end{gather*}
Let $w_0$ be a pole of $\Omega(w)$ which is $\ne0,\infty$. We associate the index $i(w_0)$
such that $w_0$ is a~zero of $x_{i(w_0)}Q_{i(w_0)}w^{\beta_{i(w_0)}}-1$ and we write
$I(w_0)=\{i_1,\dots,i_{d-1},i(w_0)\}$. Furthermore we let $\bm{w}_0=
(\zeta_1,\dots,\zeta_{i_{d-1}},w_0)$. Take the residue at $w=w_0$,
\begin{gather*}
\frac{Q_0(\zeta_1,\dots,\zeta_{d-1})w_0^{\tau_d}}{\prod_{j\not\in I(w_0)}
\big(x_jQ_j(\zeta_1,\dots,\zeta_{d-1})w_0^{\beta_j}-1\big)}\frac{1}{\beta_{i(w_0)}}.
\end{gather*}
When $w_0=\zeta_d$ we see that this differs by a factor
$\beta_{i_d}/\Delta_I=\sign(\beta_{i_d})/|\delta|$ from \eqref{res-at-point}.
Let $P$ be the set of poles $\ne0,\infty$ of $\Omega(w)$. We take the sum of the residues
of $\Omega(w)$ over all poles in $P$. We get
\begin{gather*}
\frac{1}{|\delta|}\sum_{w_0\in P}\Res_{w=w_0}\Omega(w)=\sum_{w_0\in P}\sign\big(\beta_{i(w_0)}\big)\Res_{\bm{z}=\bm{w}_0}
\omega(\bm{\tau},\bm{z}).
\end{gather*}
Without loss of generality we can assume for all $i$ that $\sign(\beta_i)>0$ if
and only if $\{i_1,\dots,i_{d-1},i\}\allowbreak\in\mathcal{I}_W$. Let $K=\{i_1,\dots,i_{d-1}\}$ and let $\mathcal{I}_K = \{I \in \mathcal{I}_W\colon K \subset I\}$ and $\mathcal{J}_K = \{I \in \mathcal{J}_W\colon K \subset I\}$.
Thus our summation becomes
\begin{gather*}
\sum_{\bm\mu\colon I(\bm\mu)\in\mathcal{I}_K}\Res_{\bm{z}=\bm{\zeta}^{\bm\mu}}\omega(\bm{\tau},\bm{z})
-\sum_{\bm{\nu}\colon I(\bm\nu)\in\mathcal{J}_K}\Res_{\bm{z}=\bm{\zeta}^{\bm\nu}}\omega(\bm{\tau},\bm{z}).
\end{gather*}
To complete our proof we need to show that $\sum_{w_0\in P}\Res_{w=w_0}\Omega(w)=0$.
Since the sum of all residues of a one variable rational function is zero, it suffices
to show that $\Res_{w=0}\Omega(w)+\Res_{w=\infty}\Omega(w)=0$. We prove that both
residues are $0$. For the residue at $w=0$ we expand~$\Omega(w)$ into a Laurent series in
$w$ times $\frac{{\rm d}w}{w}$. The support of this series is contained in the set of integers
\begin{gather*}
\ge \tau_d+\sum_{j\not\in K}\max(0,-\beta_j)=\tau_d-\sum_{j\colon \beta_j<0}\beta_j.
\end{gather*}
Since $\bm\tau$ is in $Z_B^{\circ}$ we know that there exist $\lambda_1,
\dots,\lambda_N\in(0,1)$ such that $\bm\tau=\sum_{j=1}^N\lambda_j\bm{b}_j$.
Hence we have $\tau_d=\sum_{j=1}^N\lambda_j\beta_j$ and
\begin{gather*}
\tau_d-\sum_{j\colon \beta_j<0}\beta_j =\sum_{j\colon \beta_j>0}\lambda_j\beta_j+
\sum_{j\colon \beta_j<0}(\lambda_j-1)\beta_j.
\end{gather*}
All terms in this summation are positive, hence the Laurent series expansion
of $\Omega(w)$ is in fact a~Taylor series with a zero constant term multiplied by ${\rm d}w/w$.
Hence $\Res_{w=0}\Omega(w)=0$. We~deal similarly with $w=\infty$.
\end{proof}

\begin{Lemma}\label{lem-adj-conn}
Let $\mathcal{I}$ and $\mathcal{J}$ be two cotriangulations. Then there exists a sequence of cotriangulations $\mathcal{I}_1,\dots,\mathcal{I}_N$ such that $\mathcal{I}_1 = \mathcal{I}$, $\mathcal{I}_N = \mathcal{J}$ and $\mathcal{I}_{i}$ and $\mathcal{I}_{i+1}$ are adjacent for all $i=1,\dots,N-1$.
\end{Lemma}

\begin{proof}
Let $\mathcal{I}_{\bm\rho}$ correspond to the cotriangulation with convergence direction $\bm\rho$ and $\mathcal{I}_{\bm\rho'}$ correspond to the cotriangulation with convergence direction $\bm\rho'$. Then make a continuous path $f\colon [0,1] \rightarrow \mathbb{R}^d$ such that $f(0) = \bm\rho$ and $f(1) = \bm\rho'$ which may only cross walls of the secondary fan in one point. It cannot cross lower dimensional faces of the secondary polytope. Consider the sequence $0 < t_0 < \dots < t_N < 1$ which are all points such that $f(t_i)$ is on a wall, and consider the sequence of cotriangulations
\begin{gather*}
\mathcal{I}_{f(0)},\ \mathcal{I}_{f\big(\frac{t_{0} + t_{1}}{2}\big)},\
\mathcal{I}_{f\big(\frac{t_{1} + t_{2}}{2}\big)},\ \dots,\
\mathcal{I}_{f\big(\frac{t_{N-1} + t_{N}}{2}\big)},\ \mathcal{I}_{f(1)}.
\end{gather*}
Then each consecutive cotriangulation is adjacent by definition of the path.
\end{proof}

\begin{Corollary}\label{cor-gen-tri}
Let $\mathcal{I}$ and $\mathcal{J}$ be two different cotriangulations
and suppose $\bm\tau\in Z_B^\circ$.
Then
\begin{gather*}
\sum_{\bm{\mu}\colon I(\bm{\mu})\in\mathcal{I}}\Res_{\bm{z}=\bm\zeta^{\bm\mu}}\omega(\bm{\tau},\bm{z})
=\sum_{\bm{\nu}\colon I(\bm{\nu})\in\mathcal{J}}\Res_{\bm{z}=\bm\zeta^{\bm\nu}}\omega(\bm{\tau},\bm{z}).
\end{gather*}
\end{Corollary}

\begin{proof}
Suppose $\mathcal{I}$ and $\mathcal{J}$ are adjacent cotriangulations with common wall $W$.
For the cotriangles $I \in \mathcal{I}$ such that $I \in \mathcal{J}$, there is nothing to prove as the summands on both side cancel each other out. So we are left with sums over $\bm{\mu}$ and $\bm\nu$ for corresponding cotriangles in $\mathcal{I}_W$ and $\mathcal{J}_W$ respectively. Now we simply apply Proposition \ref{lem-edge-res}.

Now suppose $\mathcal{I}$ and $\mathcal{J}$ are not adjacent cotriangulations. Then by Lemma \ref{lem-adj-conn} there exists a~sequence of adjacent cotriangulations between $\mathcal{I}$ and $\mathcal{J}$. We can now apply Proposition \ref{lem-edge-res} to each pair of adjacent contriangulations in the sequence.
\end{proof}

\section{Remarks}
\begin{Remark}
Corollary \ref{cor-gen-tri} together with Lemma~\ref{eval-res-complete} gives the final step in the proof of Theo\-rem~\ref{MB-Hermitian}
which establishes the existence of an invariant Hermitian form with respect to
$\mathcal{M}_{MB}$. The question remains whether this Hermitian form is uniquely
determined (up to a constant factor). As we know this uniqueness is equivalent to
the irreducibility of the action of $\mathcal{M}_{MB}$. In all explicit examples we
have seen so far, the Hermitian form is indeed unique.
\end{Remark}

\begin{Remark}
Recent work by Saiei Matsubara and Yoshiaki Goto \cite[Theorem 3.3]{Mat20_2}, \cite[Theorem 2.11]{Mat20} confirms the signature computation of Theorem \ref{MB-Hermitian}.
Their work does not assume
the existence of a Mellin--Barnes basis. In \cite[Theorem 3.3]{Mat20_2} they claim that the signature of the
invariant Hermitian form for any $A$-hypergeometric function with totally non-resonant
parameter vector $\bm\alpha$ with $h(\bm{\alpha})\not\in\mathbb{Z}$
is determined by the signature of
\begin{gather}
\sin\Bigg({-}\pi \sum_{i \not \in I(\bm\mu)}\gamma^{\bm\mu}_i\Bigg)\prod_{i \not \in I(\bm\mu)}\sin\big(\pi \gamma^{\bm\mu}_i\big),\qquad\bm\mu\colon\ I(\bm\mu) \in \mathcal{I},
\label{eq:Matsubara}
\end{gather}
where $\mathcal{I}$ is any cotriangulation.

In Theorem \ref{MB-Hermitian}
we see that the signature corresponds to those of
\begin{gather*}
\Delta_{I(\bm\mu)} \prod_{l \in I(\bm\mu)} (-1)^{\gamma_l^{\bm\mu}}\prod_{i \not \in I(\bm\mu)} \sin\big(\pi \gamma_i^{\bm\mu}\big),\qquad \bm\mu\colon\ I(\bm\mu) \in \mathcal{I}.
\end{gather*}
Since $\Delta_{I(\bm\mu)}>0$ we can ignore $\Delta_{I(\bm\mu)}$.
Also note that
\begin{gather*}
-\sum_{i \not \in I}\gamma^{\bm\mu}_i = -\sum_{i=1}^N \gamma_{0i} + \sum_{i \in I}\gamma^{\bm\mu}_i.
\end{gather*}
When $i \in I(\bm\mu)$ then $\gamma^{\bm\mu}_i \in \mathbb{Z}$, hence these contribute to a sign change in the leftmost $\sin$ function in \eqref{eq:Matsubara}. This sign change is exactly the product
\begin{gather*}
\prod_{l \in I(\bm\mu)} (-1)^{\gamma_l^{\bm\mu}}.
\end{gather*}
So this means we can rewrite \eqref{eq:Matsubara} to
\begin{gather*}
\sin\Bigg({-}\pi \sum_{i=1}^N\gamma_{0i}\Bigg)\prod_{l \in I(\bm\mu)} (-1)^{\gamma_l^{\bm\mu}}
\prod_{i \not \in I(\bm\mu)}\sin\big(\pi \gamma^{\bm\mu}_i\big),\qquad \bm\mu\colon\ I(\bm\mu) \in \mathcal{I}.
\end{gather*}
Note that the left-most factor equals $\sin(-\pi h(\bm\alpha))$. So when
$h(\bm\alpha)\not\in\mathbb{Z}$ we recover our result.
\end{Remark}

\subsection*{Acknowledgements}

I would like to thank Frits Beukers for his help and the anonymous referees for their comments.

\pdfbookmark[1]{References}{ref}
\LastPageEnding

\end{document}